\documentclass[11pt]{amsart}
\usepackage{amssymb, amscd, latexsym, amsmath, amsfonts, color, epsfig}

\usepackage[all]{xy}
\xyoption{dvips}

\newcommand{\R}{{\mathbb R}}

\newcommand{\Q}{{\mathbb Q}}

\newcommand{\K}{{\mathcal{K}}}
\newcommand{\W}{{\mathcal{W}}}
\newcommand{\V}{{\mathcal{V}}}

\newcommand{\Kdot}{{K^{\bullet}}}

\newcommand{\C}[1]{{C[{#1}]}}
\newcommand{\CC}[1]{{C\langle #1\rangle}}

\newcommand{\refT}[1]{Theorem~\ref{T:#1}}
\newcommand{\refC}[1]{Corollary~\ref{C:#1}}
\newcommand{\refP}[1]{Proposition~\ref{P:#1}}
\newcommand{\refD}[1]{Definition~\ref{D:#1}}

\theoremstyle{plain}
\newtheorem{thm}{Theorem}[section]
\newtheorem{prop}[thm]{Proposition}

\newtheorem{cor}[thm]{Corollary}

\theoremstyle{definition}
\newtheorem{definition}[thm]{Definition}

\theoremstyle{remark}
\newtheorem{rem}{Remark}

\begin{document}

\title{Calculus of the embedding functor and spaces of 
knots}
\author{Ismar Voli\'c}
\address{Department of Mathematics, University of Virginia, 
Charlottesville, VA}
\email{ismar@virginia.edu}
\urladdr{http://www.people.virginia.edu/\~{}iv2n}

\begin{abstract}
We give an overview of how calculus of the embedding functor can be used for the study of long knots and summarize various results connecting the calculus approach to the rational homotopy type of spaces of long knots, collapse of the Vassiliev spectral sequence, Hochschild homology of the Poisson operad, finite type knot invariants, etc.  Some open questions and conjectures of interest are given throughout.
\end{abstract}

\maketitle

{\tableofcontents}

%

%

%

%

%

%

    
\section{Introduction}\label{S:Tower}


Fix a linear inclusion $f$ of $\R$ into $\R^n$, $n\geq 3$, and let $Emb(\R, \R^n)$ and $Imm(\R, \R^n)$ be the spaces of smooth embeddings and immersions, respectively, of $\R$ in $\R^n$ which agree with $f$ outside a compact set.  It is not hard to see that $Emb(\R, \R^n)$ is equivalent to the space of based knots in the sphere $S^n$, and is known as the space of \emph{long knots}.  Let $\K^n$ be the homotopy fiber of the inclusion $Emb(\R, \R^n)\hookrightarrow Imm(\R, \R^n)$ over $f$.  Since $Imm(\R, \R^n)\simeq \Omega S^{n-1}$ (by Smale-Hirsch) and since the above inclusion is nullhomotopic \cite[Proposition 5.17]{Dev2}, we have
$$
\K^n\simeq Emb(\R, \R^n)\times \Omega^2 S^{n-1}.
$$
Spaces $\K^n$ can be thought of as spaces of long knots which come equipped with a regular homotopy to the long unknot.

The main goal of this paper is to describe and summarize the applications of the Goodwillie-Weiss calculus of the embedding functor to the study of $\K^n$.   In Section \ref{S:Tower}, we begin by introducing the Taylor tower for $\K^n$ which arises from embedding calculus and is the starting point for all other work mentioned here.  Two alternative descriptions which are better suited for computation, the mapping space model and the cosimplicial model, are given in Sections \ref{S:MappingModel} and \ref{S:CosimplicialModel}.  The connection of the cosimplicial model to the Kontsevich operad of compactified configuration spaces is the subject of Section \ref{S:KontsevichOperad}.

Some of the main results for the case $n>3$ are given in Section \ref{S:Collapses}.  Namely, the cosimplicial model for the tower gives rise to cohomology and homotopy spectral sequences, both of which collapse at the second term (Theorems \ref{T:CohomologyCollapse} and \ref{T:HomotopyCollapse}).  This completely determines the rational homotopy type of  $\K^n$, $n>3$.  Further, in Section \ref{S:Formality} we describe a combinatorial graph complex representing the $E_2$ term of the cohomology spectral sequence.  This gives a way of computing any rational cohomology group of $\K^n$, but such calculations are very hard in practice. Despite its combinatorial simplicity, the $E_2$ term is still mysterious in many ways.

To relate this to something more familiar, we recall Vassiliev's cohomology spectral sequence \cite{Vas} in Section \ref{S:Vassiliev} and describe its connection to Sinha's spectral sequence in Section \ref{S:Formality} (see \refP{SameSS} in particular).  The collapse of the latter turns out to imply the collapse of the former, settling a conjecture of Kontsevich.  A consequence is that the cohomology of $\K^n$, $n>3$, is the Hochschild homology of the degree $n-1$ Poisson operad  
(\refC{FormalityConsequences}).

Section \ref{S:Classical} is devoted to the connection between the Taylor tower and finite type invariants of $\K^3$.  Rationally, an algebraic version of the Taylor tower classifies  finite type invariants (\refT{IntroMainTheorem1}).  It is also known that the mapping space model of the Taylor tower gives the integral type two invariant (\refT{IntegralType2}) with a nice geometric interpretation by quadrisecants of the knot.  Construction of all integral finite type invariants using the Taylor tower is one of the main open questions in the subject.

Lastly, in Section \ref{S:Orthogonal} we give a brief overview of how orthogonal calculus can also be used for the study of $\K^n$.

Various open questions and conjectures are given throughout the paper.


\subsection{Acknowledgements}  The author would like to thank Dev Sinha for comments and suggestions, as well as for organizing, along with Fred Cohen and Alan Hatcher, the \emph{AIM Workshop On Moduli Spaces of Knots} for which this survey was prepared.

    
\subsection{Taylor towers for spaces of long knots arising from embedding calculus}\label{S:Tower}


One variant of calculus of functors which can be applied to $\K^n$ is  Goodwillie-Weiss \emph{embedding calculus} \cite{We, GW}.  Given a manifold $M$ and a contravariant functor $F$ from the category of open subsets of $M$ to the category of spaces or spectra, such as $Emb(M,N)$, where $N$ is another manifold, the general theory gives a \emph{Taylor tower} of fibrations
$$
F(-) \longrightarrow \big((T_{\infty}F(-)\longrightarrow \cdots \longrightarrow T_r F(-)\longrightarrow \cdots \longrightarrow T_1F(-)\big).
$$
The stage $T_rF(-)$ is the \emph{$r$th Taylor polynomial of $F$}.  In some cases, the map from $F$ to $T_{\infty}$ is an equivalence, and the tower is then said to \emph{converge}.  For the embedding functor, we have the following important result due to Goodwillie, Klein, and Weiss.

\begin{thm}[\cite{GK, GW}]\label{T:Connectivity}
If $dim(N)\!-\!dim(M)\!>\!2$, the Taylor tower for  $Emb(M,N)$ (or for the fiber of $Emb(M,N)\hookrightarrow Imm(M,N)$) converges.
\end{thm}
Note that, in the case of $\K^n$, this unfortunately says that the Taylor tower only converges for $n>3$ and not for $n=3$.  However, the tower can still be defined even for $n=3$.  

The definition of the stages $T_r$ simplifies for $\K^n$ and in fact reduces to a very concrete construction.  Namely, let 
$\{A_{i}\}$, $0\leq i\leq r$, be a collection of disjoint intervals in $\R$, and define spaces of ``punctured long knots'' as homotopy fibers
$$
\K_{S}^n=\text{hofiber} \big( Emb(\R-\bigcup_{i\in S}A_{i}, \R^{n})
\hookrightarrow
Imm(\R-\bigcup_{i\in S}A_{i}, \R^{n})\big)
$$
for each nonempty subset $S$ of $\{0, \ldots, r\}$.  Since there
are restriction maps $\K_{S}^n\to \K_{S\cup \{i\}}^n$ which commute, the
$\K_{S}^n$ can be arranged in a subcubical diagram (a cubical diagram
without the initial space) of dimension $r$.

\begin{definition}
The \emph{$r$th stage of the Taylor tower for $\K^n$}, $n\geq 3$, denoted by $T_r\K^n$, is defined to be the homotopy limit of the subcubical diagram of knots with up to $r$ punctures described above.
\end{definition}
This homotopy limit can be though of as the collection of spaces of maps of $\Delta^{|S|-1}$ into $\K_{S}^n$ fro all $S$ which are all compatible with the restriction maps in the subcubical diagram.  Section 2 of \cite{V2} has more details about these diagrams, as well as the precise definition of their homotopy limits.

Since an $r$-subcubical diagram is a face of an $(r+1)$-subcubical one, there are maps $T_{r+1}\K^n\to T_r\K^{n}$, as well as canonical maps $\K^n\to T_r\K^n$ obtained by restricting a knot to knots with up to $r+1$ punctures, which then clearly fit in the subcubical diagram.

\begin{rem}
For $r>2$, the actual limit of the subcubical diagram is equivalent to $\K^n$ itself.  One can thus think of the approach here as replacing $\K^n$, a limit of a certain diagram, by the homotopy limit of the same diagram.  This homotopy limit, although harder to define, should be easier to understand.
\end{rem}


\subsection{Mapping space model for the Taylor tower}\label{S:MappingModel}


For the following two models for the Taylor tower, the first step is to notice that punctured knots are essentially configuration spaces.  

Let $C(r)$ denote the configuration space of $r$ labeled points in $\R^n$, taken modulo the action generated by translation and scaling (since we want a compact manifold).  Then it is almost immediate (see proof of Proposition 5.13 in \cite{Dev}) that there is a homotopy equivalence
\begin{equation}\label{E:EvaluationMap}
\K_{S}^n \stackrel{\simeq}{\longrightarrow} C(|S|-1)
\end{equation}
given by evaluating the punctured knot on some points in each of the embedded arcs.  

\begin{rem}\label{R:SecondTangentialRemark}

The tangential data normally associated to an equivalence like this has been removed because our embedding spaces are fibers of inclusions of embeddings to immersions.  Otherwise the above equivalence would have to account for tangential $(n-1)$-spheres associated to the points on which the punctured knot is evaluated.  This would be the case if we worked with $Emb(\R, \R^n)$ instead of $\K^n$.  However, since most of the main results we describe
are either only known to be true or easier to state for $\K^n$ (i.e. without the tangential spheres), we have chosen to work with this space from the beginning.
\end{rem}

Because of the above equivalence, the diagrams defining the stages $T_r\K^n$ can almost be thought of as diagrams of configuration spaces.  The problem, however, is that restrictions between punctured knots have to correspond to somehow adding, or doubling, a configuration point.  To make sense of this, the configuration space first has to be compactified.  

Let $\phi_{ij}$ be the map from $C(r)$ to $S^{n-1}$ given by the 
normalized difference of $i$th and $j$th configuration point and 
$s_{ijk}$ be the map to $[0,\infty]$ given by $\vert x_{i}-x_{j}\vert/\vert 
x_{i}-x_{k}\vert.$

\begin{definition}\label{D:DevCompactifications}
Define $\C{r}$ to be the closure of the image of $C(r)$ in 
$(S^{n-1})^{{r \choose 2}}\times [0,\infty]^{r\choose 
3}$ under the product of all $\phi_{ij}\times s_{ijk}$.  Similarly 
Define 
$\CC{r}$ to be the closure of the image of $C(r)$ in 
$(S^{n-1})^{{r \choose 2}}$ under the product of all $\phi_{ij}$.
\end{definition}

These definitions were independently made in \cite{Kont, KontSoib} and in \cite{Dev1}, although Sinha was the first to explore the difference between the two spaces and to show they are both homotopy equivalent to $C(r)$ \cite[Theorem 5.10]{Dev1}.  Space $\C{r}$ is homeomorphic to the Fulton-MacPherson compactification of the configuration space \cite{AS,FM} where points are allowed to come together, while $\CC{r}$ is the quotient of this compactification by subsets of three or more points colliding along a line. 
The 
main feature of the compactifications in \refD{DevCompactifications}
is that the directions of approach of colliding configuration points are kept track of (in $\C{r}$, their relative rates of approach are also taken into account).  The stratifications of these spaces with corners have nice connections to certain categories of trees \cite{Dev1, Dev}.  An important observation is that, for configurations in $\R$, $\C{r}$ is precisely the \emph{Stasheff associahedron} $A_r$ \cite{Stash}.

The following model for the Taylor tower is convenient since it is more geometric and it interpolates between the homotopy-theoretic approach of calculus of functors and some well-known constructions such as Bott-Taubes configuration space integrals \cite{BT}.  Let $AM_r(\R^n)$ be the space of maps from $A_r$ to $\C{r}$ which are \emph{stratum-preserving} (each stratum in $A_r$ is sent to the stratum in $\C{r}$ where points with the same indices collide) and \emph{aligned} (each set of three or more points in a stratum is sent to a set of that many colinear point in a stratum) \cite[Definition 5.1]{Dev}. 

\begin{thm}\cite[Theorem 5.2]{Dev}\label{T:MappingModel}
$AM_r(\R^n)$ is homotopy equivalent to $T_r\K^n$.  Further, the map $\K^n\to AM_r(\R^n)$ given by evaluation on a knot agrees with the map $\K^n\to T_r\K^n$ in the homotopy category.  
\end{thm}

Section 5 of \cite{Dev} is devoted to the proof of this theorem.  Intuitively, the associahedron $A_r$ captures the compatibilities of the maps of simplices defining the homotopy limit $T_r\K^n$.

Using \refT{Connectivity}, we then have
\begin{cor}\label{C:KnotsMappingModel}
For $n>3$, the induced evaluation map 
$
\K^n\longrightarrow AM_{\infty}(\R^n)
$
is an equivalence.

\end{cor}


\subsection{Cosimplicial model for the Taylor tower}\label{S:CosimplicialModel}


We next want to describe a cosimplicial model for the Taylor tower, mainly because every cosimplicial space comes equipped with a cohomology and a homotopy spectral sequence which we describe in Sections \ref{S:Formality} and \ref{S:Coformality}.  This model for the tower was suggested in \cite{GKW}, but Sinha \cite{Dev} was the first to make it precise.

It turns out that spaces $\CC{m}$ admit a cosimplicial structure, while ordinary Fulton-MacPherson compactifications do not (but on the other hand, they are not manifolds with corners, which causes technical difficulties).  In general, cosimplicial diagrams are closely related to subcubical diagrams since a truncation of a cosimplicial diagram at the $r$th stage can be turned into an $r$-subcubical diagram whose homotopy limit is the $r$th partial totalization $\mbox{Tot}^r$ of the original cosimplicial diagram.  It is not true, however, that every subcubical diagram comes from a cosimplicial one, but \refT{StagesCosimplicialModel} below says that this is true in the case of the subcubical diagrams defining the Taylor tower for $\K^n$.  More about subcubical diagrams and cosimplicial spaces can be found in \cite[Section 6]{Dev}.

\begin{definition}  Let $\Kdot$ be the collection of spaces $\{ \CC{m}\}_{m=0}^{\infty}$ with doubling (coface)  maps 
$$d^i\colon \CC{m-1}\longrightarrow \CC{m}, \ \ \ 0\leq i\leq m-1,$$
which for each $i$ repeat all the vectors indexed on the $i$th configuration point, and with forgetting (codegeneracy) maps
$$s^i\colon \CC{m}\longrightarrow \CC{m-1}, \ \ \ 1\leq i\leq m,$$
which for each $i$ delete all vectors indexed on the $i$th configuration point (and relabel appropriately).
\end{definition}

\begin{thm}\label{T:StagesCosimplicialModel} For $n\geq 3$,

\begin{itemize}

\item  $\Kdot$ is a cosimplicial space \cite[Corollary 4.22]{Dev}.
\item  $\mbox{Tot}^{\, r}\Kdot\simeq T_{r}\K^n$ \cite[Theorem 1.1]{Dev2}.

\end{itemize}

\end{thm}

In analogy with \refC{KnotsMappingModel}, a consequence of this and \refT{Connectivity} is thus
\begin{cor}\label{C:KnotsCosimplicialModel}
For $n>3$, there is a homotopy equivalence $\K^n \stackrel{\simeq}{\longrightarrow} \mbox{Tot} \Kdot$ given by a collection of compatible evaluation maps.
\end{cor}

We thus have a different model for the Taylor tower given by the
sequence of spaces and fibrations

\begin{equation}\label{E:TotModel}
\K^n \longrightarrow \big((\mbox{Tot}\Kdot\longrightarrow \cdots \longrightarrow \mbox{Tot}^r\Kdot\longrightarrow \cdots \longrightarrow \mbox{Tot}^1\Kdot\big).
\end{equation}


\subsection{McClure-Smith framework and the Kontsevich operad}\label{S:KontsevichOperad}


Another way to arrive at $\Kdot$ is through the work of McClure and Smith \cite{MS} where one can associate a cosimplicial object to any operad with multiplication.  Spaces $\CC{m}$ form such an operad, called the \emph{Kontsevich operad}, which is equivalent to the little cubes operad \cite{Sal, LV2}.  Sinha then proves

\begin{thm}\cite[Theorem 1.1]{Dev2}\label{T:CosimplicialOperad}  The $r$th partial totalization of the cosimplicial space associated to the Kontsevich operad is equivalent to $T_r\K^n$.
\end{thm}

Moreover, it is in fact not hard to see that this cosimplicial space itself is equivalent to $\Kdot$.
 
 McClure and Smith also show that the totalization of the cosimplicial space arising from a multiplicative operad has an action of the little two-cubes operad \cite{MS} (see also \cite{MS2} for a more general case and a nice overview of the interplay between operads and cosimplicial spaces).  Combining this with the fact that the Taylor tower converges to $\K^n$ for $n>3$, we have

\begin{thm}\cite[Theorem 1.4]{Dev2}
For $n>3$, $\K^n$ is a two-fold loop space.
\end{thm}

\noindent
An immediate question arising from these results is:


\begin{itemize}
\item Is the two-cubes action on $\mbox{Tot}\Kdot$ compatible with the action on $\K^3$ defined by Budney \cite{B}?  Further, Budney shows that $\K^3$ is equivalent to a free little two-cubes object over the prime long knots.  Is there an equivalent freeness result for $n>3$?  In general, bringing Budney's natural and geometric two-cubes action into the picture in any way would probably be beneficial.
\end{itemize}


\section{Vassiliev spectral sequence and the Poisson operad}\label{S:Vassiliev}


We will want to relate the cohomology spectral sequence arising from $\Kdot$ to a well-known spectral sequence due to Vassiliev, who initiated the study of 
of $\K^n$ and the computation of its cohomology in \cite{Vas} by considering spaces of embeddings as complements of ``discriminants", i.e. spaces of maps with singularities.  A cohomology spectral sequence converging to $\K^n$ for $n>3$ arises from a filtration associated to the number of singularities.  Turchin gave the $E_1$ term of this spectral sequence the following convenient description.

\begin{thm}\cite[Theorem 5.11]{T}\label{T:Turchin}
The $E_1$ term of the Vassiliev cohomology spectral sequence converging to $\K^n$ is the Hochschild homology of the degree $n-1$ Poisson operad.
\end{thm} 

For the definition of the Poisson operad, see example d) in Section 1 of \cite{T} or Definition 4.10 in \cite{Dev2}, and see Section 2 of \cite{T} for the Hochschild homology of a multiplicative operad.  Another way to think of this homology is as the homology of the cosimplicial space associated to the multiplicative operad via McClure-Smith setup.

Kontsevich conjectured the collapse of Vassiliev's spectral sequence for $n\geq 3$ at $E_1$ but was only able to show this on the diagonal using what is now known as the \emph{Kontsevich integral} \cite{Kont2}.  For $n=3$, this integral gives the famous correspondence between finite type knot invariants and chord diagrams (see \refT{Kontsevich}). One of the consequences of the work described in Section \ref{S:Formality} is that Vassiliev's spectral sequence collapses at $E_1$ everywhere.


\section{Rational homotopy type of $\K^n, n>3$}\label{S:Collapses}


Following \cite{Bous, BK}, one can associate second-quadrant spectral sequences to $\Kdot$ which for $n>3$ converge to the homotopy and cohomology of $\mbox{Tot} \Kdot$ \cite[Theorems 7.1 and 7.2]{Dev} (recall that $\mbox{Tot} \Kdot\simeq\K^n$ in the same range).  
Both of these spectral sequences are now known to collapse at $E_2$ rationally.  We explain these results next and additionally in Section \ref{S:Formality} describe a graph complex giving the rational cohomology of $\K^n$.

\begin{rem}
As far as we know, the condition $n>3$ required for the convergence of the spectral sequences is independent of the same condition needed for the convergence of the Taylor tower.
\end{rem}


\subsection{Collapse of the cohomology spectral sequence}\label{S:Formality}


The cohomology spectral sequence for $\Kdot$ has
\begin{align}
E_{1}^{-p,q} & =  coker\Big(\sum (s^i)^*\colon  H^q(\CC{p-1})\longrightarrow  H^q(\CC{p})\Big), \ \ \ p,q\geq 0,  \label{E:E1} \\
d_{1} & = \sum (-1)^i (d^i)^* \colon E_{1}^{-p,q} \longrightarrow
E_{1}^{-p+1,q}. \notag
\end{align}

Sinha has shown that the $E_2$ term of this spectral sequence is also the Hochschild homology of the Poisson operad using the McClure-Smith approach mentioned in Section \ref{S:KontsevichOperad} \cite[Corollary 1.3]{Dev2}.  The main ingredient in that proof is the fact that the homology of the Kontsevich operad in $\R^n$ is the degree $n-1$ Poisson operad \cite{Coh} (see also \cite{Dev3} for an exposition of this result).  Sinha's result can also be viewed as a consequence of the following observation of Turchin, combined with \refT{Turchin}.

\begin{prop}\label{P:SameSS}
The $E_1$ term of the Vassiliev spectral sequence is isomorphic, up to regrading, to the $E_2$ term of Sinha's spectral sequence from \eqref{E:E1}.
\end{prop}

The most important ingredient in understanding this spectral sequence rationally turns out to be the {\em formality of the little cubes operad}. Recall that a space $X$ is {\em formal} if there exist rational quasi-isomorphisms of DGAs (differential graded algebras) between $(C^*(X), d\big)$ and $\big(H^*(X), 0\big)$, where $C^*$ stands
for the standard deRham-Sullivan cochain functor and $H^*(X)$ is thought of as a DGA with zero differential.  A chain complex is formal if it is quasi-isomorphic to its homology.  Formality of a diagram of spaces or an operad is defined by requiring that the quasi-isomorphisms commute with all the maps.  Kontsevich proves that the operad of chains on the little cubes operad is formal \cite{Kont} by showing that the operad of chains on the operad of Fulton-MacPherson compactifications of configuration spaces is formal.  Passing to the homotopy equivalent Kontsevich operad, and using the connection between this operad and $\Kdot$ via \refT{CosimplicialOperad}, formality ultimately gives

\begin{thm}\cite{LV}\label{T:CohomologyCollapse}
Sinha's cohomology spectral sequence whose $E_1$ term is given in \eqref{E:E1} collapses at $E_2$ for $n> 3$.
\end{thm}

The reason this is true is essentially that the vertical differential in the cohomology spectral sequence \eqref{E:E1} can be replaced by the zero differential.  Some consequences are

\begin{cor}\label{C:FormalityConsequences}  For $n>3$,
\begin{itemize}
\item
Vassiliev's spectral sequence converging to the cohomology of $\K^n$ collapses at $E_1$.
\item
The cohomology of $\K^n$ is the Hochschild homology of the degree $n-1$ Poisson operad.
\item
The DGA
$\Big(      \bigoplus\limits_{p=0}^{\infty}   s^{-p} H^*(\CC{p}) ,\,
              d_{1}=\sum_{i=0}^{p+1}(-1)^i (d^i)^*
 \Big) $, where $s^{-p}$ denotes the degree shift and $(d^i)^*$ are the maps induced by the doubling maps on cohomology,
 is a rational model for  $\K^{n}.$
 \item Sinha's cohomology spectral sequence converging to $Emb(\R,\R^n)$ collapses at $E_2$.
\end{itemize}

\end{cor}

The first statement, which proves Kontsevich's conjecture, follows from \refP{SameSS}.  The second statement follows from \refT{Turchin}, while the third, giving the DGA model for $\K^n$, is essentially the definition of the total complex of the $E_2$ page of the spectral sequence \cite{LV}.
Since long knots are $H$-spaces with addition given by ``stacking", this completely determines the rational homotopy type of $\K^n$, $n>3$.  The last statement is true because the tangential spheres, which would now have to be added to the description of the $E_1$ term (and to the DGA model for $Emb(\R,\R^n)$), contribute no differentials beyond $d_1$ as they are formal as well.

\vskip 6pt
\noindent
For computation in this spectral sequence, one can take advantage of the fact that the cohomology of configuration spaces $\CC{p}$ is well known \cite{Coh} and can be represented by certain kinds of chord diagrams.  
%
%
In more detail, let $\alpha_{ij}=\phi_{ij}^*\omega_{ij}$, $1\leq i,j\leq p$, where $\phi_{ij}$ is as before the normalized difference of points $x_i$ and $x_j$ in the configuration and $\omega_{ij}$ is the rotation-invariant unit volume form on the $(n-1)$-sphere.  Then
\begin{equation}\label{E:Model}
H^{*}(\CC{p})\cong \Lambda \alpha_{ij}/\sim,
\end{equation}
where the equivalence relations $\sim$ are
\begin{gather*}
\alpha_{ii}=0, \ \ \ 
\alpha_{ij}^{2}=0,  \ \ \ 
\alpha_{ij}=(-1)^{n}\alpha_{ji}, \ \ \ 
 \alpha_{ij}\alpha_{kl}=(-1)^{n+1}\alpha_{kl}\alpha_{ij}, \ \ \text{and} 
 \\
\alpha_{ij}\alpha_{jk}+\alpha_{jk}\alpha_{ki}+\alpha_{ki}\alpha_{ij}=0 
\ \ \text{(three-term relation)}. \label{E:Cohomology4}
\end{gather*}

An easy combinatorial way of representing this cohomology is as the vector space (over $\Q$ in the case of interest to us) generated by diagrams of $p$ labeled points, usually drawn on a line segment, with chords joining vertices $i$ and $j$ for each $\alpha_{ij}$.  An example of this correspondence is given in Figure \ref{F:CorrespondenceExample}.

\begin{figure}[h]
\begin{center}
\input{CorrespondenceExample.pstex_t}
\caption{}
\label{F:CorrespondenceExample}
\end{center}
\end{figure}

The above relations then have obvious diagram interpretations, as do the degeneracies $(s^i)^*$ and faces $(d^i)^*$.  The degeneracies are given by deleting a vertex and relabeling.  It is then not hard to see that the $E_1$ term of the spectral sequence is obtained by imposing one more relation on these chord diagrams:  \emph{Every vertex must be joined by a chord to another.}  From this, one immediately gets, for example, that the $E_1$ term has a vanishing line of slope $(1-n)/2$ \cite[Corollary 7.4]{Dev}.
(There is also an upper vanishing line studied by Turchin in \cite{T2}, so that $E_1$ is concentrated in an angle.)
The faces, which define the first differential $d_1$, are given by identifying two consecutive vertices (contracting the line segment between them).  So for example, the three generators of $E_{1}^{-4, 2(n-1)}$ and their differential are given in Figure \ref{F:Generators}.      The differential of the second diagram is zero because of the three-term relation.  Note that this is also true for the sum of the differentials of the first and the third diagram.
Since there is nothing in the image of $d_1$ in this slot because of the existence of the vanishing line, we get that, for $n>3$, $H^{2n-6}(\K^n)$ is isomorphic to $\Q^2$, generated by the second diagram (i.e. by $\alpha_{13}\alpha_{24}$) and the sum of the first and the third diagrams (i.e. by $\alpha_{12}\alpha_{34}+\alpha_{14}\alpha_{23}$) in the top line of Figure \ref{F:Generators}.  One can also easily see, for example, that both $H^{n-3}(\K^n)$ and  $H^{2n-5}(\K^n)$ are isomorphic to $\Q$, generated respectively by diagrams in Figure \ref{F:Generators3}.
  
\begin{figure}[h]
\begin{center}
\input{Generators1.pstex_t}
\vskip 12pt
\input{Generators2.pstex_t}
\caption{}
\label{F:Generators}
\end{center}
\end{figure}

\begin{figure}[h]
\begin{center}
\input{Generators3.pstex_t}
\caption{}
\label{F:Generators3}
\end{center}
\end{figure}

With this diagram combinatorics in hand, it is in principle possible to compute the $E_{2}^{-p,q}$ term for any $p,q$.  However, the computations are difficult, and the $E_2$ term is still not very well understood.  For more details about $E_1$ and $d_1$ in terms of chord diagrams, see \cite[Section 6]{V1}.

\vskip 6pt
\noindent
Some further questions are:

\begin{itemize}

\item Can we understand the combinatorics in the $E_2$ page better (with the ultimate goal of obtaining a closed form for $H^*(\K^n)$)?  
  More precisely, what is the structure of the Poisson algebra underlying this combinatorics?  What are the geometric representatives of the generators?

\item  Cattaneo,
Cotta-Ramusino, and Longoni \cite{Catt} generalize integration
techniques developed by Bott and Taubes \cite{BT, V3} in deRham theory
and produce complexes reminiscent of the rows in the spectral
sequence described here which map to the cochains on $\K^n$.  Is this map a quasi-isomorphism?
 
\item Can we say anything about torsion?



\end{itemize}


\subsection{Collapse of the homotopy spectral sequence}\label{S:Coformality}


The homotopy spectral sequence is constructed analogously, with 
\begin{align}
E^{1}_{-p,q} & = \bigcap ker\Big((s^i)_*\colon \pi_q(\CC{p+1})\longrightarrow \pi_q(\CC{p})\Big), \ \ \ p,q\geq 0,  \label{E:HE1} \\
d^{1} & =  \sum (-1)^i (d^i)_* \colon E_{1}^{-p+1,q} \longrightarrow
E_{1}^{-p,q}. \notag
\end{align}
In the rational case, homotopy groups of configuration spaces are known to form a Yang-Baxter algebra \cite{Kohno} (see also Section 2 of \cite{SS}).  This was used in \cite{SS} for some computer-aided computations, with the difficulty growing exponentially.  The authors  also show that the spectral sequence has a vanishing line, and they prove

\begin{thm}\cite[Theorem 4.7]{SS}
The Euler characteristic of the $i$th nontrivial row of the $E^1$ term is zero for $i>2$.
\end{thm}

It turns out that the notion dual to formality is useful for showing the collapse of this spectral sequence.  A space $X$ is \emph{coformal} when there exist quasi-isomorphisms connecting the rational homotopy groups of $\Omega X$ and the free Lie algebra on $C^*X$.  Coformality of a diagram or an operad is defined by as usual requiring that the quasi-isomorphisms be compatible with all the maps.  In analogy with Kontsevich's formality result, we have

\begin{thm}[\cite{ALTV}] The Kontsevich operad is coformal for $n> 3$.
\end{thm}

One then deduces a result parallel to \refT{CohomologyCollapse}.

\begin{thm}[\cite{ALTV}]\label{T:HomotopyCollapse}  Sinha's homotopy spectral sequence whose $E^1$ term is given in \eqref{E:HE1} collapses at $E^2$ for $n> 3$. 
\end{thm}

\vskip 6pt
\noindent
A further question is:

\begin{itemize}

\item Formality and coformality together say that the cohomology of long knots is essentially the free algebra on the dual of the Yang-Baxter algebra.  Understanding the combinatorial structure of this algebra should help with computations and give another point of view on the appearance of chord diagrams in the study of knots. 

\end{itemize}


\section{Case $\K^3$ and finite type invariants}\label{S:Classical}


One special case of much interest is that of classical long knots $\K^3$. Even
though one no longer has the Goodwillie-Klein-Weiss comparison of
\refT{Connectivity}, nor is it clear what the spectral sequences
converge to, the Taylor tower still provides a lot of information.

Recall that any knot invariant $V\in H^0(\K^3)$ can be extended to singular knots with $m$
transverse double points via the repeated use of the
\emph{Vassiliev skein relation}
 from Figure
\ref{Fig:Vassiliev}.
\begin{figure}[h]
\label{Fig:Vassiliev}
\begin{center}
\input{Vassiliev.pstex_t}
\caption{Vassiliev skein relation}
\end{center}
\end{figure}
Then $V$ is a \emph{type m invariant} if it vanishes identically
on knots with $m+1$ double points. Let $\V_{m}$ be the collection
of all type $m$ invariants and note that
$\V_{m-1}\!\subset\!\V_{m}$.  It turns out that $\V_{m}$ is closely related to the space generated by certain chord diagrams with $m$ chords (in fact, precisely those described in Section \ref{S:Formality} which lie on the diagonal of Sinha's spectral sequence at $E_1$ for $n=3$) modulo a certain relation (called four-term relation, which comes from $d_1$ in the spectral sequence, but also has a nice geometric interpretation).  In fact, if $\W_m$ denotes the dual of this space, then Kontsevich \cite{Kont2} (see also \cite{BN2}) proves 

\begin{thm}\label{T:Kontsevich}
$\W_m\cong\V_m / \V_{m-1}$.
\end{thm}
As mentioned in Section \ref{S:Vassiliev}, this theorem, which is the cornerstone of finite type knot theory, can be restated as the collapse of Vassiliev's spectral sequence on the diagonal at $E_1$.  
A nice selection of its various proofs can be found in \cite{BN}.  
%

The approach to finite type theory has
thus far exclusively been through integration and combinatorics.  Calculus of functors, however, might provide a topological point of view which has been missing from the theory.

The first evidence for this is the following.  After constructing the Taylor tower as described in
Section \ref{S:Tower}, one also has Sinha's cosimplicial model as well the associated Bousfield-Kan
cohomology spectral sequence from Section \ref{S:CosimplicialModel}.  However, the spectral sequence no longer necessarily converges to $\mbox{Tot}\Kdot$ but rather to an algebraic version of it, obtained by taking cochains on each of the configuration spaces $\CC{m}$ and then forming an algebraic realization of the resulting simplicial group.  (This is done essentially by collecting diagonally in the double complex whose 
vertical differential is the coboundary in $C^{*}\CC{m}$
while the alternating sum of the cofaces as usual gives the horizontal differential.)  
If one truncates the cosimplicial space at the $r$th stage, the spectral
sequence computes the cohomology of the partial algebraic totalization.  This totalization  is equivalent to the homotopy colimit of the subcubical diagram obtained by replacing the punctured knots by cochains on those spaces.  In this way, one obtains the corresponding algebraic replacement of the original Taylor tower whose stages we denote by $T^{*}_{r}\K^3$.
We then have 
\begin{thm}\cite[Theorem 6.10]{V1}\label{T:IntroMainTheorem1}
There is an isomorphism $H^{0}(T^{*}_{2m}\K^3)\cong\V_{m}.$  Moreover, all rational finite type invariants factor through the algebraic Taylor
tower for
$\K$.
\end{thm}
The factorization is through the collection of evaluation maps \eqref{E:EvaluationMap}.  There are also isomorphisms $H^{0}(T^{*}_{2m+1}\K^3)\cong H^{0}(T^{*}_{2m}\K^3)$ \cite[Section 6.2]{V1} so that all stages of the algebraic tower are accounted for.  The main ingredient in the proof of the above theorem is Bott-Taubes configuration space integrals \cite{BT, V3}.  One also recovers \refT{Kontsevich} as a corollary.


Another important connection between \emph{integral} finite type theory and embedding calculus is the work of Budney, Conant, Scannell, and Sinha \cite{BCSS}, who study the first three stages of the ordinary Taylor tower in detail.  Using the mapping space model for the tower, they further give a new geometric interpretation of the unique (up to framing) integral type 2 invariant via quadrisecants.  Their main result, which in fact holds for $Emb(\R, \R^3)$ rather than $\K^3$, is

\begin{thm}\cite[Theorem 6.3]{BCSS}\label{T:IntegralType2}
The map $\pi_0 (Emb(\R, \R^3))\longrightarrow \pi_0 (AM_3(\R^3))$ represents the unique integral additive type 2 invariant.
\end{thm}

\noindent
Some further questions are:

\begin{itemize}
\item Sinha's
cohomology spectral sequence may not converge to the
totalization of $\Kdot$ for $n=3$ because it is not clear that $C^{*}$
commutes with totalization. Further, the Taylor tower is not known
to converge to the space of knots.  Another way to say this is that the maps
\begin{equation}\label{E:SeparationMap}
H^0(\K^3)\stackrel{\alpha}{\longleftarrow}
H^0(T_\infty\K^3)=H^0\Big(\lim\limits_{r\to\infty}\underset{S\subseteq\{1,
..., r\}}{\mbox{holim}}\K_{S}^3\Big) \stackrel{\beta}{\longleftarrow}
H^0\Big(\lim\limits_{r\to\infty}\underset{S\subseteq\{1, ...,
r\}}{\mbox{hocolim}}C^* \K_{S}^3\Big)
\end{equation}
are not known to be isomorphisms.  If the spectral sequence
indeed does not converge to the desired associated graded, then
does
the genuine Taylor tower for $\K^3$ contain more information than just the finite type invariants?
The work in \cite{BCSS}
indicates that this might be the case.  Conant has further shown that type $m$ invariants show up in stage $m-1$ of the ordinary tower, rather than stage $2m$ as is the case in the algebraic tower.

\item
A very closely related question is
that of separation of knots by finite type invariants.  Since \refT{IntroMainTheorem1} states that rational finite type
invariants are all one finds on the right side of the map $\beta$,
showing that $\beta$ and $\alpha$ are surjections would settle the
question of finite type invariants separating knots.

\item As conjectured in \cite{BCSS}, can one construct all integral finite type invariants using the collection of evaluation maps and the mapping space model for the Taylor tower?  What are the analogs of quadrisecants for higher order finite type invariants?

\item Can we show collapse of Vassiliev's spectral sequence for $\K^3$ directly?  What are the convergence issues there?

\item Can one gain topological insight
into the common thread between knots in $\R^{3}$ and $\R^{n}$,
namely the Kontsevich and Bott-Taubes integral constructions of
finite type invariants? This is to be expected, since the latter
type of integrals plays a crucial role in the proof of
\refT{IntroMainTheorem1}.



\item Can we use calculus of functors to study finite type invariants of braids or 
homology spheres, for example? 
\end{itemize}


\section{Orthogonal calculus}\label{S:Orthogonal}

There is another brand of calculus of functors, \emph{orthogonal calculus}, due to Weiss \cite{We2}, which can be used for studying $\K^n$.  Here one considers covariant functors from the category of vector spaces with isometric linear inclusions to the category of spaces or spectra.  To such a functor, orthogonal calculus associates a different Taylor tower, whose layers (fibers of the maps between the stages) are spectra with actions of the orthogonal groups.

\begin{thm}[\cite{ALV}]
Let M be a smooth manifold such that $2dim(M)+1<n$.  Then the orthogonal Taylor tower for $Emb(M, \R^n)$ (or for the fiber of $Emb(M, \R^n)\hookrightarrow Imm(M, \R^n)$) splits rationally into the product of its layers.
\end{thm}
The main ingredient in the proof of this theorem is again the formality of the chains of little cubes operad combined with the interplay between embedding calculus (which views $Emb(M, \R^n)$ as a functor of $M$) and orthogonal calculus (which views $Emb(M, \R^n)$ as a functor of $\R^n$).  An immediate consequence is
\begin{cor}\label{C:OrthogonalCollapse}
The rational cohomology spectral sequence associated to the orthogonal tower for the fiber of $Emb(M, \R^n)\hookrightarrow Imm(M, \R^n)$ collapses at $E_1$ for $2dim(M)+1<n$.
\end{cor}

\begin{itemize}
\item
Combining the first statement in \refC{FormalityConsequences} with \refC{OrthogonalCollapse} for $M=\R$, we have that the $E_1$ term of the orthogonal calculus cohomology spectral sequence has to be isomorphic, up to regrading, to the $E_1$ term of the Vassiliev spectral sequence for $\K^n$, $n>3$ (although this has not yet been verified directly).  Can studying the layers of the orthogonal tower, whose explicit description is given in \cite{ALV}, give more insight into the combinatorics of the Vassiliev spectral sequence an help in computations?
\end{itemize}

%



\end{document}